\documentclass[12pt,thmsa]{article}%
\usepackage{amssymb}
\usepackage{amsfonts}
\usepackage{sw20bams}%
\usepackage{amsmath}%
\usepackage{graphicx}
\providecommand{\U}[1]{\protect\rule{.1in}{.1in}}
\begin{document}

\author{Steven R. Finch}
\title{Convex Hull of Two Orthogonal Disks}
\date{March 12, 2016}
\maketitle

\begin{abstract}
Three configurations of two perpendicular disks in $\mathbb{R}^{3}$ are
examined, the first in which the disks share centers and the other two in
which the disks touch at precisely one point. \ Volume, surface area and mean
width calculations dominate the discussion. \ Integrated mean curvature also
appears as an indirect way to compute mean width.

\end{abstract}

\footnotetext{Copyright \copyright \ 2012, 2016 by Steven R. Finch. All rights
reserved.}Our investigation begins with a theoretical question about
experimental data. \ Example 1 is the convex hull of the following two
orthogonal disks in $\mathbb{R}^{3}$:%
\[%
\begin{array}
[c]{ccccc}%
\left\{  (x,y,z):x^{2}+y^{2}\leq1\text{ \&\ }z=0\right\}   &  & \text{and} &
& \left\{  (x,y,z):x^{2}+z^{2}\leq1\text{ \&\ }y=0\right\}  .
\end{array}
\]
We can numerically evaluate the volume $VL$, surface area $AR$ and mean width
$MW$ of the corresponding solid domain in Figure 1 using \cite{Msn}:
\[%
\begin{array}
[c]{ccccc}%
VL_{1}\approx2.666, &  & AR_{1}\approx10.28, &  & MW_{1}\approx1.869.
\end{array}
\]
Example 2 is the convex hull of the two disks:%
\[%
\begin{array}
[c]{ccccc}%
\left\{  (x,y,z):x^{2}+y^{2}\leq1\text{ \&\ }z=-1\right\}   &  & \text{and} &
& \left\{  (x,y,z):x^{2}+z^{2}\leq1\text{ \&\ }y=1\right\}
\end{array}
\]
with corresponding solid domain in Figure 2 and
\[%
\begin{array}
[c]{ccccc}%
VL_{2}\approx3.141, &  & AR_{2}\approx13.92, &  & MW_{2}\approx2.277.
\end{array}
\]
Example 3 is the convex hull of the two disks:%
\[%
\begin{array}
[c]{ccccc}%
\left\{  (x,y,z):x^{2}+(y+1)^{2}\leq1\text{ \&\ }z=0\right\}   &  & \text{and}
&  & \left\{  (x,y,z):(y-1)^{2}+z^{2}\leq1\text{ \&\ }x=0\right\}
\end{array}
\]
with corresponding solid domain in Figure 3 and
\[%
\begin{array}
[c]{ccccc}%
VL_{3}\approx3.627, &  & AR_{3}\approx15.97, &  & MW_{3}\approx2.645.
\end{array}
\]
Of the nine constants, just two ($VL_{1}=8/3$ and $VL_{2}=\pi$) are readily
identifiable. \ What are exact closed-form expressions for the remaining constants?

Given $\Omega$ to be a convex body in $\mathbb{R}^{3}$, a\ \textbf{width} is
the distance between a pair of parallel $\Omega$-supporting planes. Every unit
vector $v\in\mathbb{R}^{3}$ determines a unique such pair of planes orthogonal
to $v$ and hence a width $w(v)$. Let $v$ be uniformly distributed on the unit
sphere $S^{2}\subset$ $\mathbb{R}^{3}$. Then $w$ is a random variable and its
average value is the \textbf{mean width} of $\Omega$. \ Three numerical
characteristics of $\Omega$ -- volume, surface area and mean width -- are
central to our study. These quantities, along with the Euler characteristic,
form a basis of the space of all additive continuous measures that are
invariant under rigid motions in $\mathbb{R}^{3} $.

\textquotedblleft The mean width is a new measure on three-dimensional solids
that enjoys equal rights with volume and surface area\textquotedblright%
\ \cite{Rt}, hence much of this paper is devoted to computing $MW$ for our
three examples. \ What we call the direct approach is based on the definition
of $MW$; what we call the \textit{indirect} approach utilizes a connection
between $MW$ and integrated mean curvature (often called "integral" or
\textquotedblleft total\textquotedblright\ mean curvature). \ This connection
is suggested in the materials science \cite{Hi, DH} and astrophysics
literature \cite{MBW, OHPR}; the closest claim to a proof appears in
\cite{LD}, based chiefly on \cite{MPS}. \ Our paper therefore also serves to
confirm the validity of the indirect approach for certain non-polyhedral test
cases.\footnote{On page 513 of \cite{LD}, mean curvature $\mathcal{K}$ is
defined as the \textit{average} of the two principal curvatures, but this is
inconsistent with \cite{MPS}, which takes $\mathcal{K}$ to be the
\textit{sum}. We follow \cite{MPS}, defining $2H=\mathcal{K}$. Our formula
correctly gives $MW=(\ell+\pi r)/2$ for a right circular cylinder of length
$\ell$, radius $r$ \cite{SKM, Sa}.}

\section{Example 1}

The boundary $\partial\Omega$ of the convex hull $\Omega$ here is trivially
given by the surface%
\[
z=\pm\left(  \sqrt{1-x^{2}}-\left\vert y\right\vert \right)
\]
over the planar region $x^{2}+y^{2}\leq1$. \ Let%
\[%
\begin{array}
[c]{ccc}%
\varphi(x,y)=\sqrt{1-x^{2}}-y, &  & x^{2}+y^{2}\leq1\text{ \&\ }y\geq0
\end{array}
\]
then $\varphi_{x}$, $\varphi_{y}$, $\varphi_{xx}$, $\varphi_{xy}$,
$\varphi_{yy}$ denote first/second-order partial derivatives of $\varphi$ and
\[
VL_{1}=4%
{\displaystyle\int\limits_{-1}^{1}}
{\displaystyle\int\limits_{0}^{\sqrt{1-x^{2}}}}
\varphi(x,y)\,dy\,dx=\frac{8}{3},
\]%
\begin{align*}
AR_{1}  & =%
{\displaystyle\int\limits_{\partial\Omega}}
dS=4%
{\displaystyle\int\limits_{-1}^{1}}
{\displaystyle\int\limits_{0}^{\sqrt{1-x^{2}}}}
\sqrt{1+\varphi_{x}^{2}+\varphi_{y}^{2}}\,dy\,dx=2(2+\pi)\\
& =10.2831853071795864769252867....
\end{align*}

\subsection{Indirect Approach}

Let $\partial\Omega^{+}$ denote the upper portion of $\partial\Omega$ and
$\partial\Omega^{-}$ denote the lower portion. \ On the one hand, the mean
curvature of $\partial\Omega^{+}$ is%
\[
H(x,y)=-\frac{\left(  1+\varphi_{x}^{2}\right)  \varphi_{yy}-2\varphi
_{x}\varphi_{y}\varphi_{xy}+\left(  1+\varphi_{y}^{2}\right)  \varphi_{xx}%
}{2\left(  1+\varphi_{x}^{2}+\varphi_{y}^{2}\right)  ^{3/2}}=\frac{1}{\left(
2-x^{2}\right)  ^{3/2}}
\]
over the open region $x^{2}+y^{2}<1$ \&\ $y>0$. \ It follows that%
\[
H\,dS=H(x,y)\sqrt{1+\varphi_{x}^{2}+\varphi_{y}^{2}}\,dy\,dx=\frac{1}{\left(
2-x^{2}\right)  ^{3/2}}\sqrt{\frac{2-x^{2}}{1-x^{2}}}\,dy\,dx=\frac
{dy\,dx}{\left(  2-x^{2}\right)  \sqrt{1-x^{2}}}.
\]
On the other hand, the exterior dihedral angle on the semicircular edge
$x^{2}+y^{2}=1$ \&\ $y\geq0$ is%
\[
\alpha=2\arccos\left(  \frac{1}{\sqrt{2-x^{2}}}\right)
\]
because the unit exterior normal vector to $\partial\Omega^{+}$ is%
\[
\frac{1}{\sqrt{1+\varphi_{x}^{2}+\varphi_{y}^{2}}}\left(  -\varphi
_{x},-\varphi_{y},1\right)  =\sqrt{\frac{1-x^{2}}{2-x^{2}}}\left(  \frac
{x}{\sqrt{1-x^{2}}},1,1\right)  =\frac{1}{\sqrt{2-x^{2}}}\left(  x,y,y\right)
\]
and the unit exterior normal vector to the cylinder $x^{2}+y^{2}=1$ is
$(x,y,0)$. The dot product of the two vectors is $1/\sqrt{2-x^{2}}$; we
multiply the angle by two since the dihedral angle between $\partial\Omega
^{+}$ and $\partial\Omega^{-}$ is twice the preceding angle. In terms of
arclength $s=\theta$, $0\leq\theta\leq\pi$, we have%
\[
\alpha\,ds=2\arccos\left(  \frac{1}{\sqrt{2-\cos^{2}\theta}}\right)  d\theta.
\]
The surface $\partial\Omega$ is piecewise continuously differentiable and has
$n=4$ smooth edges $\varepsilon_{j}$ with (non-constant)\ dihedral angles
$\alpha_{j}$, $1\leq j\leq n$. From the general formula
\[
MW=\frac{1}{2\pi}%
{\displaystyle\int\limits_{\partial\Omega}}
H\,dS+\frac{1}{4\pi}\sum_{j=1}^{n}%
{\displaystyle\int\limits_{\varepsilon_{j}}}
\alpha_{j}\,ds,
\]
we deduce that%
\begin{align*}
MW_{1}  & =\frac{4}{2\pi}%
{\displaystyle\int\limits_{-1}^{1}}
{\displaystyle\int\limits_{0}^{\sqrt{1-x^{2}}}}
\frac{dy\,dx}{\left(  2-x^{2}\right)  \sqrt{1-x^{2}}}+\frac{8}{4\pi}%
{\displaystyle\int\limits_{0}^{\pi}}
\arccos\left(  \frac{1}{\sqrt{2-\cos^{2}\theta}}\right)  d\theta\\
& =-\frac{1}{\pi}\left(  2\sqrt{2}\ln\left(  -1+\sqrt{2}\right)
-4\operatorname*{Li}\nolimits_{2}\left(  -1+\sqrt{2}\right)
+4\operatorname*{Li}\nolimits_{2}\left(  1-\sqrt{2}\right)  \right) \\
& =1.8697727582861870379136441...
\end{align*}
where $\operatorname*{Li}\nolimits_{2}$ is the dilogarithm
\[%
\begin{array}
[c]{ccc}%
\operatorname*{Li}\nolimits_{2}(x)=%
{\displaystyle\sum\limits_{k=1}^{\infty}}
\dfrac{x^{k}}{k^{2}}=-%
{\displaystyle\int\limits_{0}^{x}}
\dfrac{\ln(1-t)}{t}dt, &  & |x|\leq1.
\end{array}
\]

\subsection{Direct Approach}

Consider the portion of $\partial\Omega$ in the first octant only. \ In this
octant, an $\Omega$-supporting plane $P^{d}:$%
\[%
\begin{array}
[c]{ccc}%
\dfrac{a}{d}x+\dfrac{b}{d}y+\dfrac{c}{d}z=1 &  & \text{(with coefficients
}a>0,b>0,c>0\text{ and scaling factor }d>0\text{)}%
\end{array}
\]
has an associated line$\ L_{xy}^{d}:$%
\[
\dfrac{a}{d}x+\dfrac{b}{d}y=1\text{ \&\ }z=0
\]
in the $xy$-plane and an associated line $L_{xz}^{d}:$%
\[
\dfrac{a}{d}x+\dfrac{c}{d}z=1\text{ \&\ }y=0
\]
in the $xz$-plane. \ Assume WLOG\ that $a^{2}+b^{2}+c^{2}=1$. The distance of
$P^{d}$ from the origin $O$ is $d$. \ Also,%
\[%
\begin{array}
[c]{ccccc}%
\text{distance of }L_{xy}^{d}\text{ from }O\text{ is }\dfrac{d}{\sqrt
{a^{2}+b^{2}}} &  & \text{and} &  & \text{distance of }L_{xz}^{d}\text{ from
}O\text{ is }\dfrac{d}{\sqrt{a^{2}+c^{2}}}.
\end{array}
\]
The largest $d$ such that
\[%
\begin{array}
[c]{ccccc}%
L_{xy}^{d}\text{ supports }x^{2}+y^{2}=1 &  & \text{or} &  & L_{xz}^{d}\text{
supports }x^{2}+z^{2}=1
\end{array}
\]
is thus%
\[
d=\max\left\{  \sqrt{a^{2}+b^{2}},\sqrt{a^{2}+c^{2}}\right\}  .
\]
Let
\[%
\begin{array}
[c]{ccccc}%
a=\cos\theta\sin\phi, &  & b=\sin\theta\sin\phi, &  & c=\cos\phi
\end{array}
\]
where $0\leq\theta\leq\pi/2$, $0\leq\phi\leq\pi/2$. \ To ensure uniformity,
think of $(\theta,\phi)$ as possessing joint density $\frac{2}{\pi}\sin\phi$.
\ We have%
\begin{align*}
MW_{1}  & =2%
{\displaystyle\int\limits_{0}^{\pi/2}}
\,%
{\displaystyle\int\limits_{0}^{\pi/2}}
\max\left\{  \sin\phi,\sqrt{\cos^{2}\theta\sin^{2}\phi+\cos^{2}\phi}\right\}
\frac{2}{\pi}\sin\phi\,d\phi\,d\theta\\
& =\frac{8}{\pi}%
{\displaystyle\int\limits_{0}^{\pi/2}}
\,%
{\displaystyle\int\limits_{\xi(\theta)}^{\pi/2}}
\sin^{2}\phi\,d\phi\,d\theta=1.8697727582861870379136441...
\end{align*}
where%
\[
\xi(\theta)=\arccos\left(  \frac{\sin\theta}{\sqrt{2-\cos^{2}\theta}}\right)
\]
is the required solution (for $\phi$ in terms of $\theta$) of the equation%
\[
\sin\phi=\sqrt{\cos^{2}\theta\sin^{2}\phi+\cos^{2}\phi}.
\]

\section{Example 2}

The curved portions of the boundary $\partial\Omega$ of the convex hull
$\Omega$ here are given by $z=\varphi(x,y)$ and $z=\psi(x,y)$, where%
\[%
\begin{array}
[c]{ccc}%
\varphi(x,y)=\sqrt{1-x^{2}}+y-1, &  & |x|\leq1,\text{ }|y|\leq1\text{
\&\ }y\geq-\sqrt{1-x^{2}},
\end{array}
\]%
\[%
\begin{array}
[c]{ccc}%
\psi(x,y)=-\sqrt{1-x^{2}}+y-1, &  & |x|\leq1,\text{ }|y|\leq1\text{ \&\ }%
y\geq\sqrt{1-x^{2}}.
\end{array}
\]
The flat portions of $\partial\Omega$ are the two disks, one of which is given
by $z=-1$ over $x^{2}+y^{2}\leq1$. These facts contribute to the following: \
\[
VL_{2}=%
{\displaystyle\int\limits_{-1}^{1}}
{\displaystyle\int\limits_{-\sqrt{1-x^{2}}}^{\sqrt{1-x^{2}}}}
\left(  \varphi(x,y)+1\right)  \,dy\,dx+%
{\displaystyle\int\limits_{-1}^{1}}
{\displaystyle\int\limits_{\sqrt{1-x^{2}}}^{1}}
\left(  \varphi(x,y)-\psi(x,y)\right)  \,dy\,dx=\pi,
\]%
\begin{align*}
AR_{2}  & =2\pi+%
{\displaystyle\int\limits_{-1}^{1}}
{\displaystyle\int\limits_{-\sqrt{1-x^{2}}}^{1}}
\sqrt{1+\varphi_{x}^{2}+\varphi_{y}^{2}}\,dy\,dx+%
{\displaystyle\int\limits_{-1}^{1}}
{\displaystyle\int\limits_{\sqrt{1-x^{2}}}^{1}}
\sqrt{1+\psi_{x}^{2}+\psi_{y}^{2}}\,dy\,dx\\
& =2\left(  \pi+2\sqrt{2}E\left(  \tfrac{1}{2}\right)  \right) \\
& =13.9235808852350105127348109...
\end{align*}
where
\[
E(\mu)=%
{\displaystyle\int\limits_{0}^{\pi/2}}
\sqrt{1-\mu\sin(\theta)^{2}}\,d\theta=%
{\displaystyle\int\limits_{0}^{1}}
\sqrt{\dfrac{1-\mu\,t^{2}}{1-t^{2}}}\,dt
\]
is the complete elliptic integral of the second kind.

\subsection{Indirect Approach}

Let $\partial\Omega^{+}$ denote the curved portion of $\partial\Omega$
prescribed by $\varphi$ and $\partial\Omega^{-}$ denote the curved portion
prescribed by $\psi$. \ We have%
\begin{align*}
-\frac{\left(  1+\varphi_{x}^{2}\right)  \varphi_{yy}-2\varphi_{x}\varphi
_{y}\varphi_{xy}+\left(  1+\varphi_{y}^{2}\right)  \varphi_{xx}}{2\left(
1+\varphi_{x}^{2}+\varphi_{y}^{2}\right)  ^{3/2}}  & =\frac{1}{\left(
2-x^{2}\right)  ^{3/2}}\\
& =-\frac{\left(  1+\psi_{x}^{2}\right)  \psi_{yy}-2\psi_{x}\psi_{y}\psi
_{xy}+\left(  1+\psi_{y}^{2}\right)  \psi_{xx}}{2\left(  1+\psi_{x}^{2}%
+\psi_{y}^{2}\right)  ^{3/2}}%
\end{align*}
everywhere and hence%
\[
H\,dS=\frac{dy\,dx}{\left(  2-x^{2}\right)  \sqrt{1-x^{2}}}
\]
as previously. \ It follows that%
\[%
{\displaystyle\int\limits_{\partial\Omega}}
H\,dS=%
{\displaystyle\int\limits_{-1}^{1}}
{\displaystyle\int\limits_{-\sqrt{1-x^{2}}}^{1}}
\frac{dy\,dx}{\left(  2-x^{2}\right)  \sqrt{1-x^{2}}}+%
{\displaystyle\int\limits_{-1}^{1}}
{\displaystyle\int\limits_{\sqrt{1-x^{2}}}^{1}}
\frac{dy\,dx}{\left(  2-x^{2}\right)  \sqrt{1-x^{2}}}=\sqrt{2}\pi.
\]
Let $\varepsilon$ denote the circular edge $x^{2}+y^{2}=1$ \&\ $z=-1$. Clearly
$\partial\Omega^{+}\cap\varepsilon$ is the semicircle with $y=-\sqrt{1-x^{2}}$
whereas $\partial\Omega^{-}\cap\varepsilon$ is the semicircle with
$y=\sqrt{1-x^{2}}$. The exterior dihedral angle on $\varepsilon$ is%
\[
\alpha=\arccos\left(  \frac{y}{\sqrt{2-x^{2}}}\right)
\]
because the unit exterior normal vector to $\partial\Omega^{+}$,
$\partial\Omega^{-}$ is%
\[
\frac{1}{\sqrt{1+\varphi_{x}^{2}+\varphi_{y}^{2}}}\left(  -\varphi
_{x},-\varphi_{y},1\right)  =\sqrt{\frac{1-x^{2}}{2-x^{2}}}\left(  \frac
{x}{\sqrt{1-x^{2}}},-1,1\right)  =\frac{1}{\sqrt{2-x^{2}}}\left(
x,y,-y\right)  ,
\]%
\[
\frac{1}{\sqrt{1+\psi_{x}^{2}+\psi_{y}^{2}}}\left(  \psi_{x},\psi
_{y},-1\right)  =\sqrt{\frac{1-x^{2}}{2-x^{2}}}\left(  \frac{x}{\sqrt{1-x^{2}%
}},1,-1\right)  =\frac{1}{\sqrt{2-x^{2}}}\left(  x,y,-y\right)
\]
respectively and the unit exterior normal vector to the horizontal disk is
$(0,0,-1)$. The dot product of the two vectors is $y/\sqrt{2-x^{2}}$. \ An
identical argument applies for the circular edge $x^{2}+z^{2}=1$ \&\ $y=1$.
\ In terms of arclength $s=\theta$, $0\leq\theta\leq2\pi$, we obtain%
\[
2%
{\displaystyle\int\limits_{\varepsilon}}
\alpha\,ds=2%
{\displaystyle\int\limits_{0}^{2\pi}}
\arccos\left(  \frac{\sin\theta}{\sqrt{2-\cos^{2}\theta}}\right)  d\theta
=2\pi^{2}
\]
which leads to the conclusion that%
\[
MW_{2}=\frac{1}{2\pi}\left(  \sqrt{2}\pi\right)  +\frac{1}{4\pi}\left(
2\pi^{2}\right)  =\frac{1}{2}\left(  \sqrt{2}+\pi\right)
=2.2779031079814441436321660....
\]

\subsection{Direct Approach}

Consider the curved portion of $\partial\Omega$ in the halfspace $x\geq0$
only. \ In this halfspace, an $\Omega$-supporting plane $P^{d}:$%
\[%
\begin{array}
[c]{ccc}%
\dfrac{a}{d}x+\dfrac{b}{d}y+\dfrac{c}{d}z=1 &  & \text{(with coefficients
}a>0,\,b,\,c\text{ and scaling factor }d>0\text{)}%
\end{array}
\]
has associated lines$\ L_{xy}^{d}$, $L_{xz}^{d}:$%
\[%
\begin{array}
[c]{ccc}%
\dfrac{a}{d}x+\dfrac{b}{d}y=\dfrac{d+c}{d}\text{ \&\ }z=-1, &  & \dfrac{a}%
{d}x+\dfrac{c}{d}z=\dfrac{d-b}{d}\text{ \&\ }y=1.
\end{array}
\]
Assume WLOG\ that $a^{2}+b^{2}+c^{2}=1$. The distance of $P^{d}$ from the
origin $O$ is $d$. \ Also,%
\[%
\begin{array}
[c]{ccccc}%
\text{distance of }L_{xy}^{d}\text{ from }O_{z}\text{ is }\dfrac{d+c}%
{\sqrt{a^{2}+b^{2}}} &  & \text{and} &  & \text{distance of }L_{xz}^{d}\text{
from }O_{y}\text{ is }\dfrac{d-b}{\sqrt{a^{2}+c^{2}}}%
\end{array}
\]
where $O_{z}=(0,0,-1)$ and $O_{y}=(0,1,0)$. \ The largest $d$ such that one of
the unit circles is supported is thus
\[
d=\max\left\{  \sqrt{a^{2}+b^{2}}-c,\sqrt{a^{2}+c^{2}}+b\right\}  .
\]
We introduce spherical coordinates as before, but with $-\pi/2\leq\theta
\leq\pi/2$, $0\leq\phi\leq\pi$ instead. \ To ensure uniformity, think of
$(\theta,\phi)$ as possessing joint density $\frac{1}{2\pi}\sin\phi$. \ We
have%
\begin{align*}
MW_{2}  & =2\,%
{\displaystyle\int\limits_{-\pi/2}^{\pi/2}}
{\displaystyle\int\limits_{0}^{\pi}}
\max\left\{  \sin\phi-\cos\phi,\sqrt{\cos^{2}\theta\sin^{2}\phi+\cos^{2}\phi
}+\sin\theta\sin\phi\right\}  \frac{1}{2\pi}\sin\phi\,d\phi\,d\theta\\
& =\frac{2}{\pi}%
{\displaystyle\int\limits_{-\pi/2}^{\pi/2}}
\,%
{\displaystyle\int\limits_{\xi(\theta)}^{\pi}}
\left(  \sin\phi-\cos\phi\right)  \sin\phi\,d\phi\,d\theta
=2.2779031079814441436321660...
\end{align*}
where%
\[
\xi(\theta)=\pi/2+\arctan\left(  \sin\theta\right)
\]
is the required solution (for $\phi$ in terms of $\theta$) of the equation%
\[
\sin\phi-\cos\phi=\sqrt{\cos^{2}\theta\sin^{2}\phi+\cos^{2}\phi}+\sin
\theta\sin\phi.
\]

\section{Example 2 (Again)}

Vinzant, using techniques in her thesis \cite{Vzt}, computed that
$\partial\Omega$ is given implicitly by the equation
\begin{align*}
0  & =-x^{2}+2y+x^{2}y-3y^{2}+y^{3}-2z-x^{2}z+6yz\\
& +x^{2}yz-5y^{2}z+y^{3}z-3z^{2}+5yz^{2}-2y^{2}z^{2}-z^{3}+yz^{3}\\
& =\left(  y-1\right)  \left(  z+1\right)  \left(  x^{2}-2y+y^{2}%
+2z-2yz+z^{2}\right) \\
& =\left(  y-1\right)  \left(  z+1\right)  (z-\varphi(x,y))(z-\psi(x,y))
\end{align*}
verifying what we already know. \ She additionally gave an elegant parametric
representation of the curved portion in $x\geq0$:%
\[%
\begin{array}
[c]{ccccccc}%
x=\sqrt{1-v^{2}}, &  & y=1-u+uv, &  & z=-u-v+uv, &  & 0\leq u\leq1\text{
\&\ }-1\leq v\leq1
\end{array}
\]
which deserves further attention. \ In the following, we reproduce our results
from the preceding section. \ The purpose in doing so is not to torture the
reader, but rather to set the stage for Example 3 (for which a parametric
representation is the only workable method available.)\ \ The Jacobian
determinant\footnote{The fact that this has indefinite sign doesn't affect the
volume calculation.}
\[
\frac{\partial(x,y)}{\partial(u,v)}=\left\vert
\begin{array}
[c]{cc}%
x_{u} & x_{v}\\
y_{u} & y_{v}%
\end{array}
\right\vert =\frac{(-1+v)v}{\sqrt{1-v^{2}}}
\]
allows us to evaluate%
\[
VL_{2}=2%
{\displaystyle\int\limits_{0}^{1}}
{\displaystyle\int\limits_{-1}^{1}}
\left(  (-u-v+uv)+1\right)  \frac{(-1+v)v}{\sqrt{1-v^{2}}}\,dv\,du=\pi.
\]
Defining%
\[%
\begin{array}
[c]{ccc}%
E=\left(  x_{u},y_{u},z_{u}\right)  \cdot\left(  x_{u},y_{u},z_{u}\right)  , &
& G=\left(  x_{v},y_{v},z_{v}\right)  \cdot\left(  x_{v},y_{v},z_{v}\right)  ,
\end{array}
\]%
\[
F=\left(  x_{u},y_{u},z_{u}\right)  \cdot\left(  x_{v},y_{v},z_{v}\right)
\]
we have%
\begin{align*}
AR_{2}  & =%
{\displaystyle\int\limits_{\partial\Omega}}
dS=2\pi+2%
{\displaystyle\int\limits_{0}^{1}}
{\displaystyle\int\limits_{-1}^{1}}
\sqrt{EG-F^{2}}\,dv\,du\\
& =2\pi+2%
{\displaystyle\int\limits_{0}^{1}}
{\displaystyle\int\limits_{-1}^{1}}
\sqrt{\frac{(1-v)(1+v^{2})}{1+v}}\,dv\,du\\
& =2\left(  \pi+2\sqrt{2}E\left(  \tfrac{1}{2}\right)  \right)  .
\end{align*}
Defining%
\[
\mathcal{N}=\frac{\left(  x_{u},y_{u},z_{u}\right)  \times\left(  x_{v}%
,y_{v},z_{v}\right)  }{\left\vert \left(  x_{u},y_{u},z_{u}\right)
\times\left(  x_{v},y_{v},z_{v}\right)  \right\vert }=\left(  \sqrt
{\frac{1-v^{2}}{1+v^{2}}},\frac{v}{\sqrt{1+v^{2}}},-\frac{v}{\sqrt{1+v^{2}}%
}\right)  ,
\]%
\[%
\begin{array}
[c]{ccc}%
L=-\left(  x_{u},y_{u},z_{u}\right)  \cdot\mathcal{N}_{u}, &  & N=-\left(
x_{v},y_{v},z_{v}\right)  \cdot\mathcal{N}_{v},
\end{array}
\]%
\[
M=-\dfrac{1}{2}\left(  \left(  x_{u},y_{u},z_{u}\right)  \cdot\mathcal{N}%
_{v}+\left(  x_{v},y_{v},z_{v}\right)  \cdot\mathcal{N}_{u}\right)
\]
we have%
\begin{align*}%
{\displaystyle\int\limits_{\partial\Omega}}
H\,dS  & =%
{\displaystyle\int\limits_{\partial\Omega}}
\frac{EN-2FM+GL}{2\left(  EG-F^{2}\right)  }\,dS\\
& =2%
{\displaystyle\int\limits_{0}^{1}}
{\displaystyle\int\limits_{-1}^{1}}
\frac{EN-2FM+GL}{2\left(  EG-F^{2}\right)  }\sqrt{EG-F^{2}}\,dv\,du\\
& =2%
{\displaystyle\int\limits_{0}^{1}}
{\displaystyle\int\limits_{-1}^{1}}
\frac{1}{\left(  1+v^{2}\right)  ^{3/2}}\sqrt{\frac{(1-v)(1+v^{2})}{1+v}%
}\,dv\,du\\
& =2%
{\displaystyle\int\limits_{0}^{1}}
{\displaystyle\int\limits_{-1}^{1}}
\frac{1}{1+v^{2}}\sqrt{\frac{1-v}{1+v}}\,dv\,du=\sqrt{2}\pi.
\end{align*}
The semicircular edge $x^{2}+z^{2}=1$, $y=1$ \&\ $x\geq0$ corresponds to $u=0$
\&\ $-1\leq v\leq1$. \ Call this $\varepsilon$. \ The exterior dihedral angle
is%
\[
\alpha=\arccos\left(  \mathcal{N}\cdot(0,1,0)\right)  =\arccos\left(  \frac
{v}{\sqrt{1+v^{2}}}\right)
\]
and arclength $s$ satisfies%
\[
ds=\sqrt{x_{v}^{2}+z_{v}^{2}}\,dv=\frac{1}{\sqrt{1-v^{2}}}\,dv.
\]
Consequently%
\[
4%
{\displaystyle\int\limits_{\varepsilon}}
\alpha\,ds=4%
{\displaystyle\int\limits_{-1}^{1}}
\frac{1}{\sqrt{1-v^{2}}}\arccos\left(  \frac{v}{\sqrt{1+v^{2}}}\right)
dv=2\pi^{2}
\]
and%
\[
MW_{2}=\frac{1}{2\pi}\left(  \sqrt{2}\pi\right)  +\frac{1}{4\pi}\left(
2\pi^{2}\right)  =\frac{1}{2}\left(  \sqrt{2}+\pi\right)
\]
as was to be shown.

\section{Example 3}

Vinzant, using techniques in her thesis \cite{Vzt}, computed that
$\partial\Omega$ here is given implicitly by the equation%
\begin{align*}
0  & =-4x^{4}+8x^{6}-16x^{2}y+36x^{4}y-16y^{2}+40x^{2}y^{2}+15x^{4}%
y^{2}+36x^{2}y^{3}+8y^{4}\\
& +6x^{2}y^{4}-y^{6}-8x^{2}z^{2}+24x^{4}z^{2}+16yz^{2}+40y^{2}z^{2}%
-78x^{2}y^{2}z^{2}-36y^{3}z^{2}\\
& +6y^{4}z^{2}-4z^{4}+24x^{2}z^{4}-36yz^{4}+15y^{2}z^{4}+8z^{6}.
\end{align*}
One could solve for this cubic (in $z^{2}$) and proceed as earlier, laboring
against the weight of complicated expressions. \ We prefer, however, to
exploit another of her elegant parametric representations:%
\[%
\begin{array}
[c]{ccccc}%
x=u\sqrt{-v(2+v)}, &  & y=\dfrac{v(1+2uv)}{1+2v}, &  & z=\dfrac{(-1+u)\sqrt
{v(2+3v)}}{1+2v}%
\end{array}
\]
for $0\leq u\leq1$ \&\ $-2\leq v\leq-2/3$. \ The Jacobian determinant%

\[
\frac{\partial(x,y)}{\partial(u,v)}=\sqrt{\frac{-v}{2+v}}\frac{2+v+6uv+6uv^{2}%
}{(1+2v)^{2}}
\]
allows us to evaluate%
\begin{align*}
VL_{3}  & =4%
{\displaystyle\int\limits_{0}^{1}}
{\displaystyle\int\limits_{-2}^{-2/3}}
\dfrac{(-1+u)\sqrt{v(2+3v)}}{1+2v}\sqrt{\frac{-v}{2+v}}\frac{2+v+6uv+6uv^{2}%
}{(1+2v)^{2}}dv\,du\\
& =\frac{2\pi}{\sqrt{3}}=3.6275987284684357011881565....
\end{align*}
Likewise,%
\begin{align*}
AR_{3}  & =4%
{\displaystyle\int\limits_{0}^{1}}
{\displaystyle\int\limits_{-2}^{-2/3}}
\sqrt{EG-F^{2}}\,dv\,du=4%
{\displaystyle\int\limits_{0}^{1}}
{\displaystyle\int\limits_{-2}^{-2/3}}
\frac{2+v+6uv+6uv^{2}}{(1+2v)^{2}}\sqrt{-\frac{1+2v+3v^{2}}{(2+v)(2+3v)}%
}\,dv\,du\\
& =4%
{\displaystyle\int\limits_{-2}^{-2/3}}
\frac{2+4v+3v^{2}}{(1+2v)^{2}}\sqrt{-\frac{1+2v+3v^{2}}{(2+v)(2+3v)}%
}\,dv=15.9716335277272626753537899....
\end{align*}
Gosper \&\ Bickford \cite{GB} conjectured that $AR_{3}=4\sigma/\tau$, where%

\begin{align*}
\sigma & =-81i\,E\left(  \frac{1}{81}\left(  17-56i\sqrt{2}\right)  \right)
+9\sqrt{-17+56i\sqrt{2}}E\left(  \frac{1}{81}\left(  17+56i\sqrt{2}\right)
\right) \\
& +108\sqrt{2}K\left(  \frac{1}{9}\right)  -72\sqrt{2}K\left(  \frac{1}%
{81}\left(  17+56i\sqrt{2}\right)  \right) \\
& -\left(  56i+32\sqrt{2}\right)  \Pi\left(  \frac{1}{27}\left(  7-4i\sqrt
{2}\right)  ,\frac{1}{81}\left(  17-56i\sqrt{2}\right)  \right) \\
& +\left(  56i+32\sqrt{2}\right)  \Pi\left(  \frac{1}{3}\left(  7-4i\sqrt
{2}\right)  ,\frac{1}{81}\left(  17-56i\sqrt{2}\right)  \right) \\
& +72i\,\Pi\left(  \frac{1}{27}\left(  7+4i\sqrt{2}\right)  ,\frac{1}%
{81}\left(  17+56i\sqrt{2}\right)  \right) \\
& -72i\,\Pi\left(  \frac{1}{3}\left(  7+4i\sqrt{2}\right)  ,\frac{1}%
{81}\left(  17+56i\sqrt{2}\right)  \right)
\end{align*}
and $\tau=18\left(  -i+2\sqrt{2}\right)  $, where $i$ is the imaginary unit,
$E(\mu)$ was defined earlier,
\[
K(\mu)=%
{\displaystyle\int\limits_{0}^{\pi/2}}
\dfrac{1}{\sqrt{1-\mu\sin(\theta)^{2}}}\,d\theta=%
{\displaystyle\int\limits_{0}^{1}}
\dfrac{1}{\sqrt{(1-t^{2})(1-\mu\,t^{2})}}\,dt
\]
is the complete elliptic integral of the first kind and
\[
\Pi(\nu,\mu)=%
{\displaystyle\int\limits_{0}^{\pi/2}}
\dfrac{1}{\left(  1-\nu\sin(\theta)^{2}\right)  \sqrt{1-\mu\sin(\theta)^{2}}%
}\,d\theta=%
{\displaystyle\int\limits_{0}^{1}}
\dfrac{1}{(1-\nu\,t^{2})\sqrt{(1-t^{2})(1-\mu\,t^{2})}}\,dt
\]
is the complete elliptic integral of the third kind. \ A\ simplification of
$\sigma$ would be good to see someday.

\subsection{Indirect Approach}

As before,%
\[
\mathcal{N}=\frac{\left(  x_{u},y_{u},z_{u}\right)  \times\left(  x_{v}%
,y_{v},z_{v}\right)  }{\left\vert \left(  x_{u},y_{u},z_{u}\right)
\times\left(  x_{v},y_{v},z_{v}\right)  \right\vert }=\left(  \sqrt
{\frac{-v(2+v)}{1+2v+3v^{2}}},\frac{1+v}{\sqrt{1+2v+3v^{2}}},\sqrt
{\frac{v(2+3v)}{1+2v+3v^{2}}}\right)
\]
and%
\begin{align*}%
{\displaystyle\int\limits_{\partial\Omega}}
H\,dS  & =4%
{\displaystyle\int\limits_{0}^{1}}
{\displaystyle\int\limits_{-2}^{-2/3}}
\frac{EN-2FM+GL}{2\left(  EG-F^{2}\right)  }\sqrt{EG-F^{2}}\,dv\,du\\
& =-4%
{\displaystyle\int\limits_{0}^{1}}
{\displaystyle\int\limits_{-2}^{-2/3}}
\frac{3v(1+2v)^{2}}{(1+2v+3v^{2})^{3/2}(2+v+6uv+6uv^{2})}\\
& \cdot\frac{2+v+6uv+6uv^{2}}{(1+2v)^{2}}\sqrt{-\frac{1+2v+3v^{2}%
}{(2+v)(2+3v)}}\,dv\,du\\
& =-4%
{\displaystyle\int\limits_{-2}^{-2/3}}
\frac{3v}{1+2v+3v^{2}}\frac{1}{\sqrt{-(2+v)(2+3v)}}\,dv\\
& =2\sqrt{2}\pi.
\end{align*}
Let $\varepsilon$ denote the arc of the semicircle $x^{2}+(y+1)^{2}=1$, $z=0 $
\&\ $x\geq0$ that runs counterclockwise from points $(0,-2,0)$ to $\left(
2\sqrt{2}/3,-2/3,0\right)  $; this corresponds to $u=1$ \&\ $-2\leq v\leq-2/3
$. The exterior dihedral angle is%
\[
\alpha=2\arccos\left(  \mathcal{N}\cdot\left(  \sqrt{-v(2+v)},1+v,0\right)
\right)  =2\arccos\left(  \frac{1}{\sqrt{1+2v+3v^{2}}}\right)
\]
because the unit exterior normal vector to the cylinder $x^{2}+(y+1)^{2}=1$ is
$(x,y+1,0)$; the arclength $s$ satisfies%
\[
ds=\sqrt{x_{v}^{2}+y_{v}^{2}}\,dv=\frac{1}{\sqrt{-v(2+v)}}\,dv.
\]
Consequently%
\begin{align*}
4%
{\displaystyle\int\limits_{\varepsilon}}
\alpha\,ds  & =8%
{\displaystyle\int\limits_{-2}^{-2/3}}
\frac{1}{\sqrt{-v(2+v)}}\arccos\left(  \frac{1}{\sqrt{1+2v+3v^{2}}}\right)
dv\\
& =4\pi\operatorname{arcsec}(3)
\end{align*}
and therefore \
\begin{align*}
MW_{3}  & =\frac{1}{2\pi}%
{\displaystyle\int\limits_{\partial\Omega}}
H\,dS+\frac{1}{\pi}%
{\displaystyle\int\limits_{\varepsilon}}
\alpha\,ds\\
& =\sqrt{2}+\operatorname{arcsec}(3)\\
& =2.6451729797138697309366179....
\end{align*}
All nine constants exhibited (at the beginning) possess closed-form
expressions, although the result for $AR_{3}$ is partly conjectural. We had
expected that there might be required \textquotedblleft more time to develop
the languages, functions, symmetries, etc., to express the constants more
naturally\textquotedblright\ \cite{Fi0}, but this belief turned out to be
overly cautious.

\subsection{Direct Approach}

Consider the portion of $\partial\Omega$ in the first octant only. \ In this
octant, an $\Omega$-supporting plane $P^{d}:$%
\[%
\begin{array}
[c]{ccc}%
\dfrac{a}{d}x+\dfrac{b}{d}y+\dfrac{c}{d}z=1 &  & \text{(with coefficients
}a>0,\,b>0,\,c>0\text{ and scaling factor }d>0\text{)}%
\end{array}
\]
has associated lines$\ L_{xy}^{d}$, $L_{yz}^{d}:$%
\[%
\begin{array}
[c]{ccc}%
\dfrac{a}{d}x+\dfrac{b}{d}y=1\text{ \&\ }z=0, &  & \dfrac{b}{d}y+\dfrac{c}%
{d}z=1\text{ \&\ }x=0.
\end{array}
\]
Assume WLOG\ that $a^{2}+b^{2}+c^{2}=1$. The distance of $P^{d}$ from the
origin $O$ is $d$. \ Also,%
\[%
\begin{array}
[c]{ccccc}%
\text{distance of }L_{xy}^{d}\text{ from }O_{-}\text{ is }\dfrac{d+b}%
{\sqrt{a^{2}+b^{2}}} &  & \text{and} &  & \text{distance of }L_{yz}^{d}\text{
from }O_{+}\text{ is }\dfrac{d-b}{\sqrt{b^{2}+c^{2}}}%
\end{array}
\]
where $O_{-}=(0,-1,0)$ and $O_{+}=(0,1,0)$. \ The largest $d$ such that one of
the unit circles is supported is thus
\[
d=\max\left\{  \sqrt{a^{2}+b^{2}}-b,\sqrt{b^{2}+c^{2}}+b\right\}  .
\]
We introduce spherical coordinates with $0\leq\theta\leq\pi/2$, $0\leq\phi
\leq\pi/2$, obtaining
\begin{align*}
MW_{3}  & =2%
{\displaystyle\int\limits_{0}^{\pi/2}}
\,%
{\displaystyle\int\limits_{0}^{\pi/2}}
\max\left\{  \left(  1-\sin\theta\right)  \sin\phi,\sin\theta\sin\phi
+\sqrt{\sin^{2}\theta\sin^{2}\phi+\cos^{2}\phi}\right\}  \frac{2}{\pi}\sin
\phi\,d\phi\,d\theta\\
& =\frac{4}{\pi}%
{\displaystyle\int\limits_{0}^{\kappa}}
\,%
{\displaystyle\int\limits_{0}^{\xi(\theta)}}
\left(  \sin\theta\sin\phi+\sqrt{\sin^{2}\theta\sin^{2}\phi+\cos^{2}\phi
}\right)  \,\sin\phi\,d\phi\,d\theta\\
& +\frac{4}{\pi}%
{\displaystyle\int\limits_{0}^{\kappa}}
\,%
{\displaystyle\int\limits_{\xi(\theta)}^{\pi/2}}
\left(  1-\sin\theta\right)  \sin^{2}\phi\,d\phi\,d\theta\\
& +\frac{4}{\pi}%
{\displaystyle\int\limits_{\kappa}^{\pi/2}}
\,%
{\displaystyle\int\limits_{0}^{\pi/2}}
\left(  \sin\theta\sin\phi+\sqrt{\sin^{2}\theta\sin^{2}\phi+\cos^{2}\phi
}\right)  \,\sin\phi\,d\phi\,d\theta\\
& =2.6451729797138697309366179...
\end{align*}
where $\kappa=\operatorname{arccsc}(3)$ and%
\[%
\begin{array}
[c]{ccc}%
\xi(\theta)=\arccos\left(  \sqrt{\dfrac{1-4\sin\theta+3\sin^{2}\theta}%
{2-4\sin\theta+3\sin^{2}\theta}}\right)  , &  & 0\leq\theta\leq\kappa
\end{array}
\]
is the required solution (for $\phi$ in terms of $\theta$) of the equation%
\[
\left(  1-\sin\theta\right)  \sin\phi=\sin\theta\sin\phi+\sqrt{\sin^{2}%
\theta\sin^{2}\phi+\cos^{2}\phi}.
\]

\section{Related Topics}

Dirnb\"{o}ck \&\ Stachel \cite{DS} studied the convex hull of the two disks:%
\[%
\begin{array}
[c]{ccccc}%
\left\{  (x,y,z):x^{2}+\left(  y+\tfrac{\delta}{2}\right)  ^{2}\leq1\text{
\&\ }z=0\right\}  &  & \text{and} &  & \left\{  (x,y,z):\left(  y-\tfrac
{\delta}{2}\right)  ^{2}+z^{2}\leq1\text{ \&\ }x=0\right\}
\end{array}
\]
when $\delta=1$ and Ira \cite{Ira} studied the same when $\delta=\sqrt{2}$.
\ These are intermediate cases relative to our Example 3 (for which $\delta
=2$) and what is essentially Example 1 (for which $\delta=0$). \ The former
case, called an \textbf{oloid}, has volume%
\begin{align*}
VL  & =\frac{2}{3}%
{\displaystyle\int\limits_{0}^{\pi/2}}
\frac{\left(  2+\cos\theta\right)  ^{2}}{\left(  1+\cos\theta\right)
\sqrt{1+2\cos\theta}}d\theta\\
& =\frac{2}{3}\left(  -1+2\sqrt{3}E\left(  \frac{\pi}{4},\frac{4}{3}\right)
+2\sqrt{3}F\left(  \frac{\pi}{4},\frac{4}{3}\right)  \right) \\
& =3.0524184684243748566972005...
\end{align*}
and surface area $AR=4\pi$, where $E$ \& $F$ are incomplete elliptic integral
of the second \&\ first kinds respectively:%
\[
E(\phi,\mu)=%
{\displaystyle\int\limits_{0}^{\phi}}
\sqrt{1-\mu\sin(\theta)^{2}}\,d\theta=%
{\displaystyle\int\limits_{0}^{\sin\phi}}
\sqrt{\dfrac{1-\mu\,t^{2}}{1-t^{2}}}\,dt,
\]%
\[
F(\phi,\mu)=%
{\displaystyle\int\limits_{0}^{\phi}}
\dfrac{1}{\sqrt{1-\mu\sin(\theta)^{2}}}\,d\theta=%
{\displaystyle\int\limits_{0}^{\sin\phi}}
\dfrac{1}{\sqrt{(1-t^{2})(1-\mu\,t^{2})}}\,dt
\]
(of course, $F(\pi/2,\cdot)=K(\cdot)$). \ A feature of the oloid is that each
of the two interlocking orthogonal circles $C_{\text{left}}$ \&
$C_{\text{right}}$ intersects the center of the other. \ Define constants
$\omega=\arcsin\left(  \sqrt{\sqrt{2}-1}\right)  $ and%
\begin{align*}
\gamma & =\frac{1}{2\sqrt{2}}%
{\displaystyle\int\limits_{0}^{\pi/2}}
\frac{\left(  2+\sqrt{2}\cos\theta\right)  ^{2}\left(  1+\sqrt{2}\cos
\theta+\cos^{2}\theta\right)  }{\left(  1+\sqrt{2}\cos\theta\right)  ^{2}%
\sqrt{1+2\sqrt{2}\cos\theta+\cos^{2}\theta}}d\theta\\
& =\frac{\sqrt{2}-1}{2}+E(\omega,-1)+\sqrt{2}\left(  \Pi(-\sqrt{2}%
-1,\omega,-1)+\Pi(\sqrt{2}-1,\omega,-1)-F(\omega,-1)\right)
\end{align*}
where $\Pi$ is the incomplete elliptic integral of the third kind:
\[
\Pi(\nu,\phi,\mu)=%
{\displaystyle\int\limits_{0}^{\phi}}
\dfrac{1}{\left(  1-\nu\sin(\theta)^{2}\right)  \sqrt{1-\mu\sin(\theta)^{2}}%
}\,d\theta=%
{\displaystyle\int\limits_{0}^{\sin\phi}}
\dfrac{1}{(1-\nu\,t^{2})\sqrt{(1-t^{2})(1-\mu\,t^{2})}}\,dt.
\]
The latter case, called a \textbf{two-circle roller}, has volume
\[
VL=\frac{8}{3\sqrt{2}}\gamma=3.2818194874496894190321933...
\]
and surface area%
\[
AR=8\gamma=13.9235808852350105127348109....
\]
As the name \textquotedblleft roller\textquotedblright\ suggests, the authors
of \cite{DS, Ira} were raising issues of a physical/mechanical nature. \ We
suspect that an algebraic approach based on \cite{Vzt} could supplant some of
their technical arguments leading to $VL$ \&\ $AR$. \ Exact expressions for
$MW$ remain open, as far as is known.

A\ picture in \cite{CR} depicts the circles $C_{\text{left}}$ \&
$C_{\text{right}}$ but not within our context of convex hulls (taut rubber
sheets spanning wire frames). \ The context instead is about \textit{minimal
surfaces} (elastic soap films spanning the same), for which mean curvature is
zero everywhere. \ Many questions suggest themselves.

\section{Acknowledgements}

Wouter Meeussen's package ConvexHull3D.m was helpful to me in preparing this
paper \cite{Msn}. \ He kindly extended the software functionality at my
request. Cynthia Vinzant generously computed the parametric representations
for both Examples 2 and 3 based on \cite{Vzt}. I am also thankful to R.
William Gosper, Neil Bickford, Roland Roth, Rolf Schneider, Hiroshi Ira, Thinh
Le, Qiang Du, Frank Sottile and Tina Mai for their helpful correspondence.

\section{Addendum}

Axel Vogt improved upon the conjecture in \cite{GB} and successfully obtained
\[
AR_{3}=\frac{4}{3}\left[  9E\left(  \frac{1}{9}\right)  -8K\left(  \frac{1}%
{9}\right)  +8\Pi\left(  -\frac{1}{3},\frac{1}{9}\right)  \right]
\]
via a change of variables and Maple. \ It is possible to reduce the
complicated expression $4\sigma/\tau$ further:%
\begin{align*}
\sigma & =-243iE\left(  \frac{1}{81}\left(  17-56i\sqrt{2}\right)  \right)
+27\sqrt{-17+56i\sqrt{2}}E\left(  \frac{1}{81}\left(  17+56i\sqrt{2}\right)
\right) \\
& +162\sqrt{2}K\left(  \frac{1}{9}\right)  -162\sqrt{2}K\left(  \frac{1}%
{81}\left(  17+56i\sqrt{2}\right)  \right) \\
& +6\left(  8i-7\sqrt{2}\right)  \Pi\left(  \frac{1}{27}\left(  7-4i\sqrt
{2}\right)  ,\frac{1}{81}\left(  17-56i\sqrt{2}\right)  \right) \\
& +2\left(  104i+71\sqrt{2}\right)  \Pi\left(  \frac{1}{3}\left(  7-4i\sqrt
{2}\right)  ,\frac{1}{81}\left(  17-56i\sqrt{2}\right)  \right)  ,
\end{align*}%
\[
\tau=54\left(  -i+2\sqrt{2}\right)
\]
using numerics and Mathematica, but Bill Gosper's question (on symbolic
transformations between elliptic integrals to link these) remains unanswered.

\begin{tabular}
[c]{lll}
& Steven R. Finch & \\
& Dept. of Statistics & \\
& Harvard University & \\
& Cambridge, MA, USA & \\
& \textit{steven\_finch@harvard.edu} &
\end{tabular}
%

\begin{figure}[ptb]%
\centering
\includegraphics[
height=3.1324in,
width=6.1583in
]%
{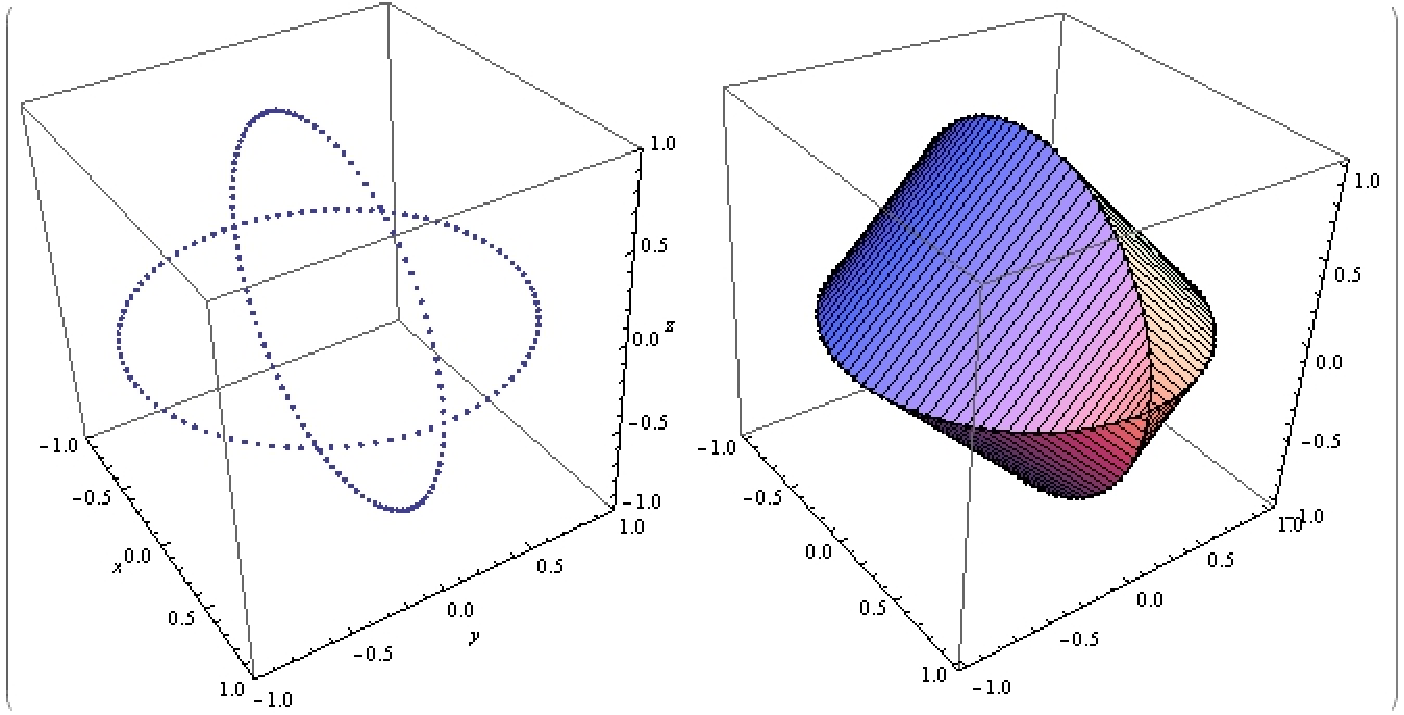}%
\caption{Centers coincide}%
\end{figure}
%

\begin{figure}[ptb]%
\centering
\includegraphics[
height=2.9793in,
width=6.1324in
]%
{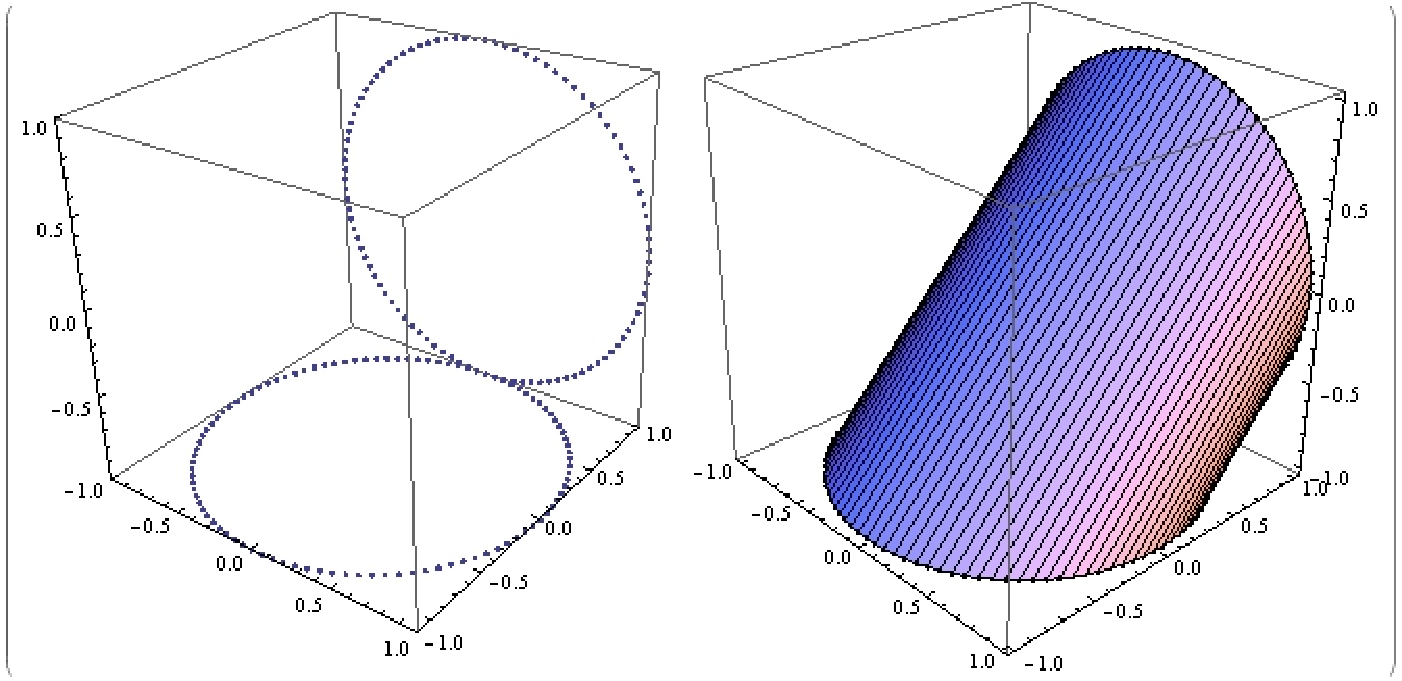}%
\caption{Centers on diagonal}%
\end{figure}
%

\begin{figure}[ptb]%
\centering
\includegraphics[
height=6.7403in,
width=4.6328in
]%
{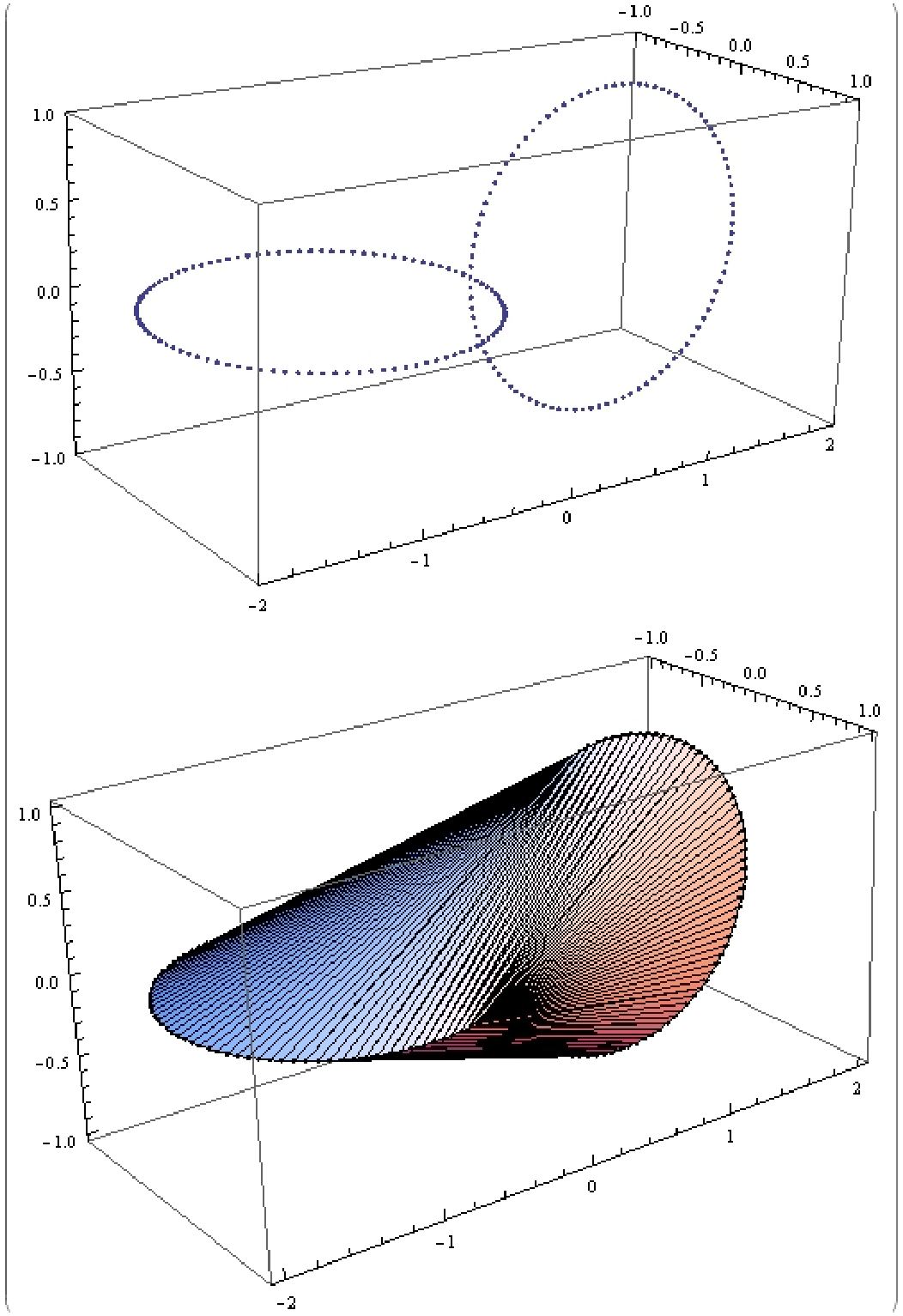}%
\caption{Centers on bisecting line}%
\end{figure}

\end{document}